\documentclass[fleqn]{mat01}
\usepackage{times,mathtimy,amssymb,latexsym}
\begin{document}

\setcounter{page}{507} \firstpage{507}

\renewcommand\theequation{\arabic{section}.\arabic{equation}}

\def\d{\mbox{\rm d}}

\newtheorem{theore}{Theorem}
\renewcommand\thetheore{\arabic{section}.\arabic{theore}}
\newtheorem{theor}[theore]{\bf Theorem}
\newtheorem{lem}[theore]{\it Lemma}
\newtheorem{propo}[theore]{\rm PROPOSITION}
\newtheorem{coro}[theore]{\rm COROLLARY}
\newtheorem{definit}[theore]{\rm DEFINITION}
\newtheorem{probl}[theore]{\it Problem}
\newtheorem{exampl}[theore]{\it Example}
\newtheorem{popp}[theore]{\it Proof of Proposition}

\def\pomp{\trivlist \item[\hskip \labelsep{\it Proof of Malliavin's theorem.}]}
\def\rema{\trivlist \item[\hskip \labelsep{\it Remark.}]}
\def\remas{\trivlist \item[\hskip \labelsep{\it Remarks.}]}

\title{Malliavin calculus of Bismut type without probability}

\markboth{R\'{e}mi L\'{e}andre}{Malliavin calculus of Bismut type
without probability}

\author{R\'{E}MI L\'{E}ANDRE}

\address{Institut de Math\'ematiques, Universit\'e de Bourgogne,
21000 Dijon, France\\
\noindent R.C.T.P. Central Visayan Institute, Jagna, Bohol, Philippines\\
\noindent E-mail: Remi.leandre@u-bourgogne.fr\\[1.2pc]
\noindent {\it Dedicated to Professor \vspace{-1pc}Sinha}}

\volume{116}

\mon{November}

\parts{4}

\pubyear{2006}

\Date{}

\begin{abstract}
We translate in semigroup theory Bismut's way of the Malliavin
calculus.
\end{abstract}


\maketitle

\section{Introduction}

There are many types of infinite dimensional analysis (see works
of Hida, Fomin, Albeverio, Elworthy, Berezanskii, etc) but one of
the specificities of Malliavin calculus is that it can be applied
to diffusions. Namely, one of the specificity of Malliavin
calculus is to complete the classical differential operations on
the Wiener space, such that functionals which belong to all the
Sobolev spaces of Malliavin calculus are in general only almost
surely defined, because there is no Sobolev imbedding theorem in
infinite dimension. Diffusions, almost surely defined, belong to
all the Sobolev spaces of Malliavin\break calculus.

By using this functional analysis approach, Malliavin
\cite{Ma$_{1}$} got a probabilistic proof of Hoermander's theorem
(see \cite{IW,Me,St,Ma$_{1}$,Ma$_{2}$,Nu} for a pedagogical
introduction).

Bismut has avoided the heavy apparatus of functional analysis of
Malliavin calculus, in order to prove again Hoermander's theorem
by probabilistic methods\break \cite{B$_{1}$,No}.

We show that Bismut's mechanism can be suitably interpreted in
terms of semigroup theory. We avoid using probability theory, but
this work is the translation in semigroup theory of the work of
Bismut. Remarkable formulas can be seen with the intuition of
probability. Let us remark that it is not the first case that
probability can explain easily magical formulas: the mysterious
rescaling of Getzler's proof of the Index theorem can be easily
interpreted in terms of probability theory in the work of
L\'eandre \cite{L$_{3}$}. We refer to the survey of L\'eandre
\cite{L$_{2}$} for various probabilistic proofs of the Index
theorem, including Bismut proof \cite{B$_{2}$}.

For the sake of simplicity, we work in the elliptic case because
in this case, the method of L\'eandre \cite{L$_{1}$,L$_{4}$} of
inversion of the Malliavin matrix can be easily interpreted in
terms of semigroup theory.

We use the classical results of differentiability of the solutions
of parabolic equations which depend on a parameter, which are
coming in stochastic analysis from the stochastic flow theorem.

\section{Cameron--Martin--Girsanov--Maruyama formula in semigroup theory}

Let us consider some vector fields $X_i$, $i=0,\dots,m$ on $R^d$
with bounded derivatives of all order. Let $L$ be the Hoermander
type operator
\begin{equation}
Lf = X_0 + 1/2 \sum_{i>0}X_i^2
\end{equation}
acting on smooth bounded functions on $R^d$. It can be written as
\begin{equation}
Lf = \langle X_0, Df\rangle  + 1/2 \sum_{i>0}\langle DX_iX_i,
Df\rangle  + 1/2 \sum_{i>0}\langle  X_i, D^2f, X_i\rangle.
\end{equation}
In (2.1), vector fields are considered as first order differential
operators and in (2.2) vectors fields are considered as smooth
applications from $R^d$ into $R^d$. Let us consider the generator
\begin{align}
L^h  = L + \sum_{i>0}h_t^iX_i,
\end{align}
where $t \rightarrow h_t^i$ are smooth bounded functions which do
not depend on $x$. $L^h$ generates an inhomogeneous Markov
semigroup $P^h$ acting on bounded continuous functions on $R^d$.

Let us consider on $R^{d+1}$ some vector fields
\begin{equation}
\tilde{X}_i^t = (X_i, h_t^iu)
\end{equation}
and the generator acting on smooth functions of $\tilde{f}$ on
$R^{d+1}$:
\begin{align}
\tilde{L}^h(\tilde{f}) &= \langle  X_0, \tilde{D}\tilde{f}\rangle
+ 1/2 \sum_{i>0}\langle DX_iX_i, \tilde{D}\tilde{f}\rangle\nonumber\\[.4pc]
&\quad\, + 1/2 \sum_{i>0}\langle \tilde{X}_i^t,
\tilde{D}^2\tilde{f}, \tilde{X}_i^t\rangle.
\end{align}
It generates a semigroup $\tilde{P}^h$ operating on the bounded
continuous functions on $R^{d+1}$. In the sequel, for the
integrability conditions, we refer to the Appendix.

\begin{theor}[\! (Quasi-invariance)]
\begin{equation}
P_t^hf(x) = \tilde{P}_t^h[uf](x,1).
\end{equation}
\end{theor}

\begin{proof}
Since the vector fields $\tilde{X}_i^t$ are linear in $u$, we have
\begin{equation}
\tilde{P}^h[uf](x,u_0) =\tilde{P}^h[uf](x,1)u_0
\end{equation}
for any bounded continuous $f$ on $R^d$ such that
\begin{equation}
\tilde{L}^h\tilde{P}^h[uf](\cdot,\cdot)\vert(x,1) =
L^h\tilde{P}^h[uf](\cdot,1)\vert (x).
\end{equation}
Therefore the result arises by using uniqueness of the solution of
the parabolic equation associated to $L^h$. \hfill  $\diamondsuit$
\end{proof}

\begin{rema}
In order to see from where this remarkable formula comes, we use
the stochastic analysis. Let $w_t^i$ be $m$ Brownian  motions and
let $Y = X_0+ 1/2\sum_{i>0}DX_iX_i$ the vector field on $R^d$. We
consider the following stochastic differential equation in It\^o
sense on $R^d$ starting from $x$:
\begin{equation}
\delta x_t^h(x) = Y(x_t^h(x))\d t + \sum_{i>0}h_t^iX_i(x_t^h(x))\d
t + \sum_{i>0}X_i(x_t^h(x))\delta w_t^i.
\end{equation}
It is a classical result of stochastic analysis that:
\begin{equation}
P_t^hf(x) = E[f(x_t^h(x))].
\end{equation}
On the other hand, let us consider the It\^o equation on $R^{d+1}$
starting from $(x,1)$:
\begin{equation}
\delta\tilde{x}_t(x) = Y(\tilde{x}_t(x)\d t + \sum
\tilde{X}_i^t(\tilde{x}_t(x))\delta w_t^i.
\end{equation}
We have
\begin{equation}
\tilde{P}^h_t[F](x,1) = E[F(\tilde{x}_t(x))]
\end{equation}
if $F$ is a bounded continuous function on $R^{d+1}$. Moreover,
the classical Girsanov formula shows us that, by using stochastic
calculus we have
\begin{equation}
E[f(x_t^h(x))] = E[F(\tilde{x}_t(x))],
\end{equation}
where $F(x,u) = uf(x)$.

Let  us consider the vector field $\bar{X}_i^h = (X_i, h_t^i)$ and
the generator on $R^{d+1}$ acting on smooth functions $\tilde{f}$
on $R^{d+1}$:
\begin{equation}
\bar{L}^h\tilde{f} = \langle X_0, \tilde{D}f\rangle + 1/2
\sum_{i>1}\langle DX_iX_i, \tilde{D}\tilde{f}\rangle  + 1/2
\sum_{i>0}\langle \bar{X}_i^h, \tilde{D}^2\tilde{f},
\bar{X}_i^h\rangle.
\end{equation}

It generates a semigroup $\bar{P}^h$ acting on the bounded
continuous functions on $R^{d+1}$.
\end{rema}

\begin{theor}[\! (Elementary integration by parts formula)]
\begin{equation}
\int_0^tP_{t-s}\sum_{i>0}h_u^iX_i[P_sf]\d s  =
\bar{P}_t^h[uf](x,0).
\end{equation}
\end{theor}

\begin{proof}
We have $\bar{P}^h_t[uf](x,u_0) = A_t(x) u_0 + B_t(x)$ because
$\partial/\partial u$ commute with $\bar{L}^h$. Therefore,
\begin{equation}
\bar{P}^h_t[uf](x,u_0) = P_t[f](x)u_0 + \bar{P}^h_t[uf](x,0)
\end{equation}
such that
\begin{equation}
\frac{\partial}{\partial t}\bar{P}^h[uf](\cdot,\cdot)\vert (x,0) =
L\bar{P}^h[uf](\cdot,0)\vert (x) + \sum_{i\geq 1}h_t^i\langle X_i,
P_t[f](x)\rangle
\end{equation}
with starting condition 0.

On the other hand, \hbox{$F(t,x) =
\int_0^tP_{t-s}[\sum_{i>0}h_s^iX_i[P_s[f]]]\d s$} is a solution of
the parabolic equation:
\begin{equation}
\frac{\partial}{\partial t}F(t,x) = L F(t,x) + \sum_{i>0}h_t^iX_iP_t[f](x)
\end{equation}
with starting condition 0. The result arises by unicity of the
solution of this parabolic equation. \hfill $\diamondsuit$
\end{proof}

\begin{rema}
Let us show from where this formula comes by stochastic analysis.
In (2.9), we put a small $\lambda$ before $h_t^i$ and we get two
processes $x_t^\lambda(x)$ and $\tilde{x}_t^\lambda(x)$. We get
the formula:
\begin{equation}
E\left[Df(x_t^0)\frac{\partial}{\partial \lambda}x_t^0(x)\right] =
E\left[DF(\tilde{x}_t^0) \frac{\partial}{\partial
\lambda}\tilde{x}_t^0(x)\right].
\end{equation}
We recognize in this last expression the quantity
\begin{equation}
f(x_t^0(x)) \int_0^t\sum_{i>0}h_t^i\delta w_t^i
\end{equation}
which gives the second term of (2.15). On the other hand,
$\frac{\partial}{\partial \lambda}x_t^0(x)$ is a solution of the
stochastic differential Stratonovitch equation:
\begin{align}
\d \frac{\partial}{\partial \lambda}x_t(x) &=
DX_0(x_t^0)\frac{\partial}{\partial \lambda}x_t(x)\d t + \sum
DX_i(x_t(x))\frac{\partial}{\partial \lambda}x_t(x) \d w_t^i\nonumber\\[.4pc]
&\quad\,+ \sum h_t^iX_i(x_t^0(x))\d t
\end{align}
starting from 0. It can be solved by the method of variation of
constant. Let $U_t$ be the solution of the matricial equation
starting from $I$:
\begin{equation}
\d U_t = DX_0(x_t^0(x))U_t\d t + \sum DX_i(x_t^0(x))U_t\d w_t^i.
\end{equation}
It is a classical result of stochastic analysis that:
\begin{equation}
\frac{\partial}{\partial \lambda}x_t^0(x) = U_t
\int_0^tU_s^{-1}h_s^iX_i(x_s^0(x))\d s.
\end{equation}
But (see Lemma~3.2 below), we have
\begin{equation}
\sum_{i>0}h_t^iX_iP_t[f](x) =\sum_{i>0}
E[Df(x_t^0(x))U_t](x,I)X_i(x).
\end{equation}
The result follows by doing the change of variable $s \rightarrow
t-s$.
\end{rema}

\section{Malliavin's theorem in semigroup theory}

For the integrability condition, we refer to the appendix. Let us
consider the vector fields on $R^d \times Gl(R^d) \times M^d =
V^d$ where $Gl(R^d)$ is the space of invertible matrices on $R^d$,
and $M^d$ the space of matrices on $R^d$: \setcounter{equation}{0}
\begin{equation}
\hat{X}_i = (X_i, DX_i U , 0)
\end{equation}
and
\begin{equation}
\hat{X} = \left(0,0, \sum_{i>0} \langle
U^{-1}X_i,\cdot\rangle^2\right).
\end{equation}

Let us consider the semigroup operating on continuous bounded
functionals on $V^d$ generated by $\hat{L}$:
\begin{align}
\hat{L}\hat{f} &= 1/2 \sum_{i>0} \langle \hat{X}_i,
\hat{D}^2\hat{f}, \hat{X_i}\rangle\nonumber\\[.4pc]
&\quad\,+ 1/2 \sum_{i>0}\langle
D\hat{X}_i\hat{X}_i,\hat{D}\hat{f}\rangle + \langle \hat{X}_0,
\hat{D}\hat{f}\rangle  + \langle \hat{X}, \hat{D}\hat{f}\rangle.
\end{align}
We consider $\hat{P}_t$ the semigroup associated to $\hat{L}$. It
is associated to the system of stochastic differential equation
(in Stratonovitch form):
\begin{align}
\d x_t(x) &= X_0(x_t(x))\d t +\sum_{i>0}X_i(x_t(x))\d w_t^i\nonumber\\[.4pc]
\d U_t(x) &= DX_0(x_t(x))U_t\d t + \sum_{i>0}Dx_i(x_t(x))U_t\d w_t^i\nonumber\\[.4pc]
\d V_t &= \sum_{i>0}\langle U_t^{-1}X_i (x_t(x),\cdot\rangle ^2
\end{align}
starting from $(x,U,V)$. If we consider $U = Id$ and $V=0$, $U_t$
corresponds to $Dx_t(x)$ and $V_t$ to the so-called Malliavin
matrix.

\setcounter{theore}{0}
\begin{theor}[\! {[Ma$_{\bf 1}$,\,Ma$_{\bf 2}$]}]
If $\hat{P}_t[V^{-p}](x,I,0) < \infty$ for all $p${\rm ,} $P_tf(x)
= \int_{R^d}p_t(x,y)$ $f(y)\d y$ where $y \rightarrow p_t(x,y)$ is
smooth positive.
\end{theor}

Let us study $X_iP_t[f]$.

\begin{lem}
Let $\bar{P}_t$ be the semigroup acting on bounded continuous
functions on $R^d \times Gl(R^d)$ associated to $\bar{L} =
\bar{X}_0 + 1/2 \sum \bar{X}_i^2$ where $\bar{X}_i = (X_i, DX_i
U)$. We get
\begin{equation}
DP_t[f](x) = \bar{P}_t[Df U](x,I).
\end{equation}
\end{lem}

\begin{proof}
We write
\begin{equation}
\psi_t(x,U_0) = \bar{P}_t[Df U](x, U_0) = \bar{P}_t[DfUU_0](x,I).
\end{equation}
In $U_0 = I$, we can compute $\bar{L}\psi_t$. We get
\begin{align}
&1/2\sum_{i>0}\langle \bar{X}_i,v\bar{D}^2\psi_t,\bar{X}_i\rangle
+ 1/2\sum_{i>0}\langle \overline{DX}_i\bar{X}_i,
\bar{D}\psi_t\rangle + \langle \bar{X}_0,
\bar{D}\psi_t\rangle\nonumber\\[.4pc]
&\quad\,= 1/2 \sum_{i>0}\langle X_i, D^2\psi_t, X_i\rangle  + 1/2
\sum_{i>0}\langle D X_iX_i, D\psi_t\rangle\nonumber\\[.4pc]
&\qquad\,   +\sum_{i>0}\langle X_i, D\psi_t,DX_i\rangle + 1/2
\sum\langle DX_iDX_i,\psi_t\rangle  + \sum_{i>0}\langle
D^2X_iX_i,\psi_t\rangle\nonumber\\[.4pc]
&\qquad\,+ \langle X_0, D\psi_t\rangle + \langle DX_0,
\psi_t\rangle
\end{align}
such that
\begin{align}
\frac{\partial}{\partial t} \psi_t(x,I) &= 1/2\sum \langle X_i,
D^2\psi_t, X_i\rangle   + 1/2 \sum \langle DX_iX_i,
D\psi_t\rangle\nonumber\\[.4pc]
&\quad\, + 1/2\sum \langle  X_i, D\psi_t, DX_i\rangle  + 1/2 \sum
\langle DX_iDX_i, \psi_t\rangle \nonumber\\[.4pc]
&\quad\, + 1/2\sum \langle D^2X_iX_i, \psi_t\rangle + \langle X_0,
D\psi_t\rangle  + \langle DX_0, \psi_t\rangle
\end{align}
with initial condition $Df$.

Moreover, $P_tf$ satisfies the equation
\begin{align}
\frac{\partial}{\partial t}P_tf = 1/2 \sum_{i>0}\langle X_i,
D^2P_tf, X_i\rangle  + 1/2\sum_{i>0}\langle DX_iX_i,
DP_tf\rangle\!+\!\langle X_0, DP_tf\rangle
\end{align}
such that $\frac{\partial}{\partial t}DP_tf$ satisfies (3.7) with
the same initial condition $Df$. Therefore we get the result by
the unicity of the solution of (3.7). \hfill $\diamondsuit$
\end{proof}

\begin{rema}
Let us show from where this remarkable formula comes. We have
\begin{equation}
P_tf(x), = E[f(x_t(x))]
\end{equation}
such that
\begin{equation}
DP_tf(x) = E[Df(x_t(x))Dx_t(x)].
\end{equation}
This result arises by the considerations following (3.4).

Let us proceed as in \cite{No}. We define on $R^{d_1}\times\cdots
\times R^{d_k}$ some vector fields
\begin{equation}
X_i^{\rm tot} = (X_i^1(x_1),\dots,
X_i^j(x^1,\dots,x^j),\dots,X_i^k(x^1,\dots,x^k)),
\end{equation}
where
\begin{align}
X_i^k(x^1,\dots,x^k) &=
X_{1,i}^k(x^1,\dots,x^{k-1})x^k\frac{\partial } {\partial x^k} +
X^k_{2,i}(x^1,\dots,x^k)\nonumber\\[.3pc]
&\quad\,+ X_{3,i}^k(x^1,\dots,x^{k-1}),
\end{align}
where $X_{1,i}^k$ has bounded derivatives of all orders,
$X_{2,i}^k= X_{2,i}^k(x^1,\dots,x^{k-1})\frac{\partial}{\partial
x_k}$ has  derivatives of all orders with polynomial growth and
$X_{3,i}^k$ has derivatives with polynomial growths.

We can define a semigroup $P^{{\rm tot},k}$ associated to $1/2
\sum_{i>0}(X_i^{{\rm tot},k})^2 + X_0^{{\rm tot},k}$. We get if
$k' < k$,
\begin{equation}
P^{{\rm tot},k}[f^{{\rm tot}, k'}](x^{{\rm tot},k}) = P^{{\rm
tot}, k'}[f^{{\rm tot},k'}](x^{{\rm tot},k'}).
\end{equation}
$x^k$ has to be seen as a matrix if $X_{1,i}^k $ is not equal to
zero.

Equation~(3.14) can be seen by using stochastic analysis because
$P^{{\rm tot},k}$ is associated to a step by step system of
stochastic differential equations. Moreover, for all $p$ we get as
follows.
\end{rema}

\setcounter{theore}{2}
\begin{propo}$\left.\right.$\vspace{.5pc}

\noindent {\it $P^{{\rm tot},k}(\vert x^k\vert^p]< \infty$ and if
$X_{2,i}^k(x^1,\dots,x^k) = 0 = X_{3,i}^k(x^1,\dots,x^{k-1})${\rm
,}
\begin{equation}
P^{{\rm tot},k}(\vert(x^k)^{-1}\vert^p](x^1,\dots,x^{k-1},I) <
\infty.
\end{equation}}
\end{propo}

We refer to the Appendix for the proof of this proposition.

Instead of considering the generator $X_0 +1/2 \sum_{i>0}X_i^2$ by
$\sum X_i h^i_t$ where $h_t^i$ is deterministic, we consider the
perturbation by $\langle\phi(x),h\rangle^iX_i$ where $\phi$ is
smooth bounded with derivatives of polynomial growths.

\looseness -1 We get a semigroup $P_t^\lambda$ associated to the
generator $X_0+1/2\!\sum X_i^2 +
\lambda\!\sum_{i>0}\langle\phi(x)\!, h\rangle^i X_i$.

\begin{lem}
Let $\bar{X}_i = (X_i, DX_iU)${\rm ,} $i = 0, 1,\dots,m$ and
$\tilde{X}_0 = (0, \sum X_i\langle \phi(x), h_t^i\rangle)$ and
$\bar{P}'_t$ be the semigroup associated to
$1/2\sum_{i>0}\bar{X}_i^2 + \bar{X}_0 + \tilde{X}_0$. We get
\begin{equation}
\frac{\partial }{\partial \lambda}P_t^0[f](x) = \bar{P}'_t[Df
U](x,0).
\end{equation}
\end{lem}

\begin{proof}
The integrability conditions are satisfied by Proposition~3.4.

Let $\bar{P}_t$ be the semigroup associated to $1/2 \sum
\bar{X}_i^2 + \bar{X}_0$. If the Volterra expansion converges for
the $C^k$ uniform norm on each compact, we get
\begin{align}
&\bar{P}'_t[DfU](x,0)\nonumber\\[.4pc]
&\quad\,=\sum(-1)^n\int_{0<s_1<\cdots<s_n<t}\bar{P}_{s_1}\tilde{X_0}
\dots\tilde{X}_0\bar{P}_{t-s_n}[Df U](x,0)\d s_1\dots\d s_n.
\end{align}
But $\bar{P}_t[DfU](x, U_0)$ is linear in $U_0$ and
$\tilde{X}_0U_0= \sum X_i\langle \phi(x), h_t\rangle^i = Y_i$
which does not depend on $U_0$. We deduce that
\begin{equation}
\bar{P}'_t[Df U](x,0) = -\int_0^tP_{s_1}\sum_{i>0}Y_i
\bar{P}_{t-s_1}[DfU](x,I)
\end{equation}
which is the formula of Lemma~3.2. \hfill $\diamondsuit$
\end{proof}

\begin{rema}
Let us show from where this formula comes. $P_t^\lambda$ is
associated to the stochastic differential equation in the
Stratonovitch sense:
\begin{equation}
\d x_t^\lambda = X_0(x_t^\lambda)\d t +
\sum_{i>0}X_i(x_t^\lambda)\d w_t^i + \lambda
\sum_{i>0}X_i(x_t^\lambda)\langle \phi(x_t^\lambda),h_t^i\rangle\d
t
\end{equation}
such that
\begin{equation}
\frac{\partial}{\partial \lambda}P_t^0[f](x) =
E\left[Df(x_t^0(s)\frac{\partial } {\partial
\lambda}x_t^0(x)\right].
\end{equation}
But $\frac{\partial}{\partial \lambda}x_t^0$ satisfies the
stochastic differential equation in the Stratonovitch sense
starting from 0:
\begin{align}
\d\frac{\partial}{\partial \lambda}x_t^0 &=
DX_0(x_t^0)\frac{\partial}{\partial \lambda}x_t^0 +
\sum_{i>0}DX_i(x_t^0)\frac{\partial}{\partial \lambda} x_t^0\d
w_t^i\nonumber\\[.4pc]
&\quad\,+ \sum_{i>0}X_i(x_t^0)\langle \phi(x_t^0),h_t^i\rangle\d t
\end{align}
and the couple of $\big(x_t^0, \frac{\partial}{\partial
\lambda}x_t^0\big)$ is associated to $\bar{P}'_t$.\pagebreak
\end{rema}

\begin{propo}\hskip -.3pc {\rm \cite{B$_{1}$}}$\left.\right.$\vspace{.5pc}
{\it
\begin{align}
\bar{P}'_t[DfU](x,0) = Q_t[fu](x,0),
\end{align}
where $Q_t$ is the semigroup generated by
\begin{equation}
\tilde{L}\tilde{f} = 1/2 \sum_{i>0}\langle\tilde{X}_i,
\tilde{D}^2\tilde{f}, \tilde{X}_i\rangle + 1/2\sum_{i>0}\langle
DX_iX_i, \tilde{D}\tilde{f}\rangle
 + \langle X_0, \tilde{D}\tilde{f}\rangle
\end{equation}
on $R^{d+1}$ where
\begin{equation}
\tilde{X}_i = (X_i, \langle\phi(x), h\rangle^i).
\end{equation}}
\end{propo}

We get the following.

\begin{proof}
We remark that the vector fields involved with $\tilde{L}$ commute
with ${\partial}/{\partial u}$ such that $Q_t[fu](x, u_0) = A(x)
u_0 + B(x)$. We remark therefore that
\begin{equation}
Q_t[fu](x,u_0) = Q_t[fu](x,0) + P_t[f](x)u_0
\end{equation}
such that
\begin{equation}
\frac{\partial}{\partial t}Q_t[fu](x,0) = LQ_t[fu](x,0) + \sum
Y_iP_t[f](x).
\end{equation}
Therefore the result is as in the proof of Theorem~2.2. \hfill
$\diamondsuit$
\end{proof}

In the previous formula, $\bar{P}_t'[DfU](x,0)$ is a scalar. We
would like to get a vector. In Lemma~3.4, we choose $\tilde{X_0} =
X_i^t(U^{-1}X_i)$ where $U$ is chosen according to Lemma~3.2. We
get with this extension
\begin{equation}
\bar{P}''_t[Df{V}](x,0) = Q_t[fu](x,I,0),
\end{equation}
where $u$ is a vector in $Q_t[fu](x,I,0)$ (see Lemma~3.7 for the
definition of $\bar{P}''$).

\begin{lem}
If $\bar{P}''_t[\vert  V^{-1}\vert^p](x,I,0)$ is finite for all
$p${\rm ,} $f \rightarrow P_t[f](x)$ has a smooth density.
\end{lem}

\begin{proof}
We can use Proposition~3.5 to the extended semigroup of
Proposition~3.3 where we replace $\phi(x)$ by $^t(U^{-1}X_i)$. We
apply Proposition~3.5 to $(V,f) \rightarrow f{V}^{-1}$. We get
\begin{equation}
\bar{P}''_t[D(f{V}^{-1})V](x,0) = P_t[Df](x) + \hbox{terms}.
\end{equation}
We iterate this procedure. We get if $\bar{P''}_t[\vert
{V}^{-1}\vert^p](x,I,0) < \infty$ for all $p$ that
\begin{equation}
P_t[D^rf](x) \leq C_r\Vert f \Vert_\infty
\end{equation}
for all $r$ for the supremum norm $\Vert \cdot \Vert_\infty$ on
functions on $R^d$. Therefore the result. \hfill $\diamondsuit$
\end{proof}

\begin{pomp}
We have if $\hat{P}_t[\vert V^{-1}\vert^p](x,I,0) < \infty$ for
all $p$ that $\bar{P}''_t[\vert {V}^{-1}\vert^p](x,I,0) < \infty$
for all $p$. For that, we use the following lemma.
\end{pomp}

\begin{lem}
Let $\tilde{X}_0'', = (0,0,Y )$ where $Y$ depends on the previous
variables and has derivatives with polynomial growth. Let
$\tilde{X}_i = (X_i, DX_iU, DX_iV)$. Let $\bar{P}''_t $ be the
semi-group associated with $1/2\sum \tilde{X}_i^2 + \tilde{X}_0 +
\tilde{X}''_0$. We have if $f$ is a homogeneous polynomial in
$V${\rm ,} $\bar{P}''_t[f](x,I, V_0) = \hat{P}_t[f(UV+UV_0)](x, I,
0)$ where we have replaced for $\hat{P}_t$ $\langle
U^{-1}X_i,\cdot\rangle^2$ by $U^{-1}Y $ in {\rm (3.2)}.
\end{lem}

\begin{proof}
Let $\tilde{P}_t$ be the semigroup associated to $1/2 \sum
\tilde{X}_i^2 + \tilde{X}_0$. It transforms a homogeneous
polynomial in $V$ into a homogeneous polynomial in $V$ of same
order (see (3.32)). Therefore, we have
\begin{align}
&\bar{P}''_t[f](x,I, V_0)\nonumber\\[.4pc]
&\quad\,= \sum (-1)^n
\int_{0<s_1<\cdots<s_n<t}\tilde{P}_{s_1}Y\tilde{P}_{s_2-s_1}Y\dots\tilde{P}_{t-s_n}[f](x,I,
V_0).
\end{align}

We put
\begin{equation}
\int_{0<s_1<\cdots<s_n<t}\tilde{P}_{s_1}Y\tilde{P}_{s_2-s_1}Y\dots\tilde{P}_{t-s_n}[f](x,I,
V_0) =I_n(s_1,\dots,s_n).
\end{equation}
Let us recall that $\hat{P}_t$ is a Markov semigroup. Let
$\hat{E}^{s_1,\dots,s_n,1}$ be the law of $x_{s_1},
U_{s_1},\dots,x_{s_n}$, $U_{s_n}, x_1,U_1$, starting from $(x,I)$
according this semigroup. We recognize in $I_n(s_1,\dots,s_n)$,
\begin{equation}
\hskip -4pc
\hat{E}^{s_1,\dots,s_n,1}[D^nf(U_1V_0)U_1U_{s_1}^{-1}Y(x_{s-1},U_{s_1})U_1U_{s_2}^{-1}Y(x_{s_2},U_{s_2})\cdots
U_1U_{s_n}^{-1}Y(x_{s_n},U_{s_n})].
\end{equation}
Therefore the series (3.30) is finite and (3.30) is valid.
But this last expression is nothing else but
\begin{equation}
\bar{P}_{s_1}\hat{X}_0\bar{P}_{s_2-s_1}\cdots\hat{X}_0
\bar{P}_{t-s_n}[f(\cdot,\cdot(V+V_0))](x, I, 0).
\end{equation}
Therefore the series (3.30) is equal to
\begin{align}
&\sum(-1)^n\int_{0<s_1<\cdots<s_n<t}\bar{P}_{s_1}\hat{X}_0
\bar{P}_{s_2-s_1}\cdots\hat{X}_0
\bar{P}_{t-s_n}[f(\cdot,\cdot(V+V_0))](x, I, 0)\nonumber\\[.4pc]
&\quad\,= \hat{P}_t[f(\cdot,\cdot(V+V_0))](x, I, 0).
\end{align}
Therefore the result. \hfill $\diamondsuit$
\end{proof}

\begin{rema}
Let us show from where this formula comes. $\bar{P}''_t$  is
associated to the system of stochastic Stratonovitch differential
equation,
\begin{align}
\d x_t(x) &= X_0(x_t(x) \d t + \sum_{i>0}X_i(x_t)\d w_t^i,\nonumber\\[.4pc]
\d U_t &= DX_0(x_t(x))U_t\d t + \sum_{i>0}DX_i(x_t(x))U_t\d
w_t^i,\nonumber\\[.4pc]
\d V_t &= DX_0(x_t(x))V_t\d t + \sum_{i>0}DX_i(x_t(x))V_t\d w_t^i
+ Y(x_t(x),U_t)\d t
\end{align}
starting from $(x, I, V_0)$. We can solve the last equation by the
method of the variation of constant, and we find that
\begin{equation}
V_t = U_t\left(V_0+\int_0^tU_s^{-1}Y(x_s(x),U_s)\d s\right).
\end{equation}
Therefore the result.
\end{rema}

\section{Inversion of the Malliavin matrix in semigroup theory}

In Theorem~3.1, $V$ is called Malliavin's covariance quadratic
form. To simplify this work, we will do the following elliptic
hypothesis assumption in $x$:
\setcounter{equation}{0}
\begin{equation}
\sum_{i >1}\langle X_i(x), \ \xi\rangle^2 \geq C \vert \xi \vert^2
\end{equation}
for some $C>0$.

\setcounter{theore}{0}
\begin{lem}
If $\vert \xi \vert = 1$,
\begin{equation}
\hat{P}_t(\vert V \xi\vert < Ct)(x, I, 0) < C < 1.
\end{equation}
\end{lem}

\begin{proof}
We introduce a function $g$ strictly decreasing, convex, from
$[0,\infty[$ into $[0,1]$ equals to 1 in 0 and tending to 0 at
infinity. We consider the function $F(\cdot)\hbox{:}\ s
\rightarrow \hat{P}_s\big[g\big(\frac{\vert V \xi\vert }{
t}\big)\big](x,I,0)$. It has a derivative in 0 in $-C/t$ and a
second derivative bounded by $C/t^2$. Moreover $F(0) = 0$ and
$F(Ct)<C<1$ for some $t$. This shows the result. \hfill
$\diamondsuit$
\end{proof}

\begin{propo}$\left.\right.$\vspace{.5pc}

\noindent {\it $\hat{P}_t[\vert V \xi\vert < \epsilon](x,I,0) \leq
C_p \epsilon^p$ for all $p$ uniformly in $\vert \xi\vert = 1$.
{\rm (}We say in such a case that $\hat{P}_t[\vert V \xi \vert
<\epsilon](x,I,0) = o(\epsilon^\infty)$.{\rm )}}
\end{propo}

\begin{proof}
We get for a big $C$,
\begin{equation}
\hat{P}_t(\vert U^{-1}\vert > C](x,I,0) = o(t^\infty).
\end{equation}
In order to show that, we choose a positive function $g = 0$ in a
neighborhood of $I$ and equal to $1$ far from $I$. We have, by
using the parabolic equation satisfied by $\hat{P}_t$
\begin{equation}
\frac{\partial^r}{\partial t^r}\bar{P}_0[g(U^{-1})](x,I) = 0
\end{equation}
for all $r$. Therefore the result.

The same result holds for $P_t(\vert\cdot-x\vert > C)(x)$.
Moreover,
\begin{equation}
\hat{P}_t[\vert V \xi\vert < \epsilon](x,I,0) \leq
\hat{P}_{\epsilon^\alpha}[\vert V \xi\vert \leq \epsilon](x,I,0)
\end{equation}
for $\alpha \leq 1$.

We slice $[0, \epsilon^\alpha]$ in $\epsilon^{-\beta}$ intervals
with $\alpha + \beta < 1$. By the previous lemma
\begin{equation}
\sup_{\vert y-x\vert <C, \vert U ^{-1}\vert <
C}P_{\epsilon^{\alpha +\beta}}[\vert V \xi\vert \leq
\epsilon](y,U,0)<C<1.
\end{equation}

We deduce by Markov property that
\begin{equation}
\hat{P}_{\epsilon^\alpha}[\vert V \xi\vert](x,I,0) <
C^{\epsilon^{-\beta}} = o(\epsilon^\infty).
\end{equation}
\hfill $\diamondsuit$
\end{proof}

\begin{theor}[\!]
$\hat{P}_t[\vert V \vert^{-p}](x,I,0) < \infty$ if $t > 0$.
\end{theor}

\begin{proof}
We remark that
\begin{equation}
\hat{P}_t(\vert V \vert^p](x,I,0) < \infty
\end{equation}
for all $p>0$ (Proposition~3.3).

We choose $\epsilon^{-\beta}$ points $\xi_i$ on the sphere of
$R^d$ such that
\begin{align}
&\hat{P}_t[\vert V^{-1}\vert > \epsilon](x,I,0)\nonumber\\[.4pc]
&\quad\,\leq \sum \hat{P}[\vert V \xi_i\vert < \epsilon](x,I,0) +
\hat{P}_t[\vert V \vert < \epsilon^{-\gamma}](x,I,0) =
o(\epsilon^\infty)
\end{align}
for some suitable $\gamma > 0$. \hfill $\diamondsuit$
\end{proof}

As a corollary, we get the following.

\begin{theor}[\!]
Under {\rm (4.1),} $f \rightarrow P_t[f](x)$ has a smooth density
$p_t(x,y)$.
\end{theor}

\section*{Appendix}

\setcounter{section}{3} \setcounter{theore}{2}
\begin{popp}{\rm
We work by induction on $k$. We choose for $k = 1$ a smooth
function $g(u) = 0$ in $u= 1$ and equal to $\vert u \vert$ when
$u$ goes to infinity, with bounded derivatives. We get for two
constant independent on $C$ that

\renewcommand\theequation{A.\arabic{equation}}

\setcounter{equation}{0}
\begin{align}
&\left\vert \frac{\partial}{\partial t}P^{{\rm tot},1}[\vert u
\vert^{2r}\exp[-g(u)/C]](x)\right\vert\nonumber\\[.4pc]
&\quad\,\leq K_1 + K_2P^{{\rm tot},1}[\vert u
\vert^{2r}\exp[-g(u)/C]](x)
\end{align}
for $r\in N$. The result follows by using the Gronwall lemma when
$C \rightarrow \infty$.

We split up the equation giving $P[f(x_k)]$ in a equation giving a
linear matrix and an equation which depends only on the previous
terms as in Lemma~3.7. We get $P_1^{{\rm tot},k}[f(x_k)](\cdot,0)
= \hat{P}_t[f(\hat{u}_k\hat{v}_k)](\cdot,I,0)$ as in Lemma~3.7.
(We start from 0 in order to simplify the exposition.) We get
$\hat{P}^{{\rm tot},k}[\vert \hat{u}_k\vert^{-p}](\cdot,I) <
\infty$ because
\begin{align}
&\left\vert\frac{\partial}{\partial t}\hat{P}^{{\rm tot},k}[\vert
\hat{u}_k\vert^{-2r}\exp[-\hat{g}(\hat{u}_k)/C]](I)\right\vert\nonumber\\[.4pc]
&\quad\,\leq K_1 + K_2\hat{P}^{{\rm tot},k}[\vert
\hat{u}_k\vert^{-2r}\exp[-\hat{g}[\hat{u}_k]/C]](I),
\end{align}
where $\hat{g}(\hat{u}_k)$ is a function with bounded derivatives
at infinity, bounded at infinity and equal to $\vert
\hat{u}_k\vert^{-1}$ near 0. The result arises by Gronwall lemma
and making $C \rightarrow \infty$.

Moreover, if we do the change of variable $u_k \rightarrow
u_k^{-1}$, we still get a semigroup governed by a generator of the
same type. This shows that
\begin{equation}
\hat{P}^{{\rm tot},k}[\vert \hat{u}_k^{-1}\vert^{-2r}](\cdot,I) <
\infty.
\end{equation}
On the other hand,
\begin{equation}
\frac{\partial}{\partial t}\hat{P}^{{\rm to}t,k}[\vert
\hat{u}_k^{-1}\vert^{2r}\bar{g}(\hat{u}_k)/C)](\cdot,I)\vert \leq
K_1 + K_2\hat{P}^{{\rm tot},k}[\vert
\hat{u}_k^{-1}\vert^{2r}\bar{g}(\hat{u}_k)/C)](\cdot,I),
\end{equation}
where $\bar{g}(\hat{u}_k)$ is a smooth function with values in
[0,1] equal to 1 in 0 and equal in a neighborhood of infinity to
$\vert \hat{u}_k^{-1}\vert^{-K}$ for a big $K$.

We use Gronwall lemma and by making $C \rightarrow \infty$ we
deduce that $\hat{P}^{{\rm tot},k}[\vert
\hat{u}_k^{-1}\vert^{2r}](\cdot,I) < \infty$ and symmetrically
that $\hat{P}^{{\rm tot},k}[\vert \hat{u}_k\vert^{2r}](\cdot,I) <
\infty$.

In order to estimate $\hat{P}^{{\rm tot},k}[\vert v_k\vert^{2r}]$
we proceed in a similar but simpler way.} \hfill $\diamondsuit$
\end{popp}

\section*{Acknowledgement}

The author would like to thank the Research Center of Theoretical
Physics, Central Visayan Institute of Jagna, Bohol, Philippines
for its kind hospitality where this work was done.

\end{document}